\documentclass[11pt, reqno]{amsart}
\usepackage{amsfonts,latexsym}
\usepackage{amsmath}
\usepackage{amscd}
\usepackage{float,amsmath,amssymb,mathrsfs,bm,multirow,graphics}
\usepackage[dvips]{graphicx}
\usepackage[percent]{overpic}
\usepackage{amsaddr}
\usepackage[numbers,sort&compress]{natbib}

\newcommand{\nequation}{\setcounter{equation}{0}}

\newcommand{\R}{{\Bbb R}}

\newcommand{\C}{{\Bbb C}}
\newcommand{\Z}{{\Bbb Z}}

\newcommand{\re}{\mathrm{Re\,}}
\newcommand{\im}{\mathrm{Im\,}}

\newtheorem{theorem}{Theorem}[section]
\newtheorem{proposition}[theorem]{Proposition}

\newtheorem{remark}[theorem]{Remark}

\newtheorem{figuretext}{Figure}

\input epsf
\title[Linearizable boundary value problems]
{Linearizable boundary value problems for the elliptic sine-Gordon and the elliptic Ernst equations}

\author{J. Lenells}
\address{Department of Mathematics, KTH Royal Institute of Technology, \\ 100 44 Stockholm, Sweden.}
\email{jlenells@kth.se}

\author{A. S. Fokas}
\address{Department of Applied Mathematics and Theoretical Physics, \\ University of Cambridge, Cambridge CB3 0WA, United Kingdom.}
\email{t.fokas@damtp.cam.ac.uk}

\begin{document}

\begin{abstract} 
\noindent
By employing a novel generalization of the inverse scattering transform method known as the unified transform or Fokas method, it can be shown that the solution of certain physically significant boundary value problems for the elliptic sine-Gordon equation, as well as for the elliptic version of the Ernst equation, can be expressed in terms of the solution of appropriate $2 \times 2$-matrix Riemann--Hilbert (RH) problems. These RH problems are defined in terms of certain functions, called spectral functions, which involve the given boundary conditions, but also unknown boundary values. For arbitrary boundary conditions, the determination of these unknown boundary values requires the analysis of a nonlinear Fredholm integral equation. However, there exist particular boundary conditions, called linearizable, for which it is possible to bypass this nonlinear step and to characterize the spectral functions directly in terms of the given boundary conditions. Here, we review the implementation of this effective procedure for the following linearizable boundary value problems: (a) the elliptic sine-Gordon equation in a semi-strip with zero Dirichlet boundary values on the unbounded sides and with constant Dirichlet boundary value on the bounded side; (b) the elliptic Ernst equation with boundary conditions corresponding to a uniformly rotating disk of dust; (c) the elliptic Ernst equation with boundary conditions corresponding to a disk rotating uniformly around a central black hole; (d) the elliptic Ernst equation with vanishing Neumann boundary values on a rotating disk.
\end{abstract}

\maketitle



\section{Introduction} \nequation
A novel method for analyzing initial-boundary value problems for integrable nonlinear evolution PDEs in one spatial dimension was introduced in 1997 \cite{F1997} (see also \cite{F2002}) and was further developed by several authors (see for example \cite{FIS2005, BFS2006, FI2004, BS}). This method, known as the unified transform or Fokas method, is based on ideas of the inverse scattering transform and it expresses the solution $q(x,t)$ of a given initial-boundary value problem in terms of the solution of a matrix Riemann--Hilbert (RH) problem. This problem involves an explicit $(x,t)$-dependence in the form $\exp [ikx-i\omega (k) t]$, where $\omega (k)$ is the dispersion relation of the associated linearized evolution PDE, as well as certain functions of the spectral parameter $k \in \mathbb{C}$, called spectral functions. The main difficulty with initial-boundary value problems, as opposed to initial value problems, is that whereas in the latter case the spectral functions are defined in terms of the given initial conditions, in the former case, the spectral functions, in addition to the given initial and boundary conditions, also involves \emph{unknown} boundary values. These unknown functions can be characterized in terms of the given data via the so-called \emph{global relation}. Significant progress in the analysis of the global relation was achieved in \cite{BFS2003} and \cite{F2005}, where it was shown that the unknown boundary values can be expressed explicitly in terms of the given data and certain eigenfunctions; however, these eigenfunctions \emph{do} depend on the unknown boundary values, thus the problem of expressing the unknown boundary values in terms of the given data remains nonlinear. 

It was shown in \cite{F2002} that for a particular class of boundary conditions, called linearizable, it is possible to bypass the above nonlinear step and express the spectral functions explicitly in terms of the given data; linearizable boundary value problems (BVPs) are discussed in \cite{BFS2004, Cau2015, C2018, FLkdv, LFgnls, F2004}.

The new method was implemented for integrable \emph{nonlinear elliptic} PDEs in two dimensions in \cite{FLP2013} (see also \cite{P, PP}) and \cite{LFernst, L, LP2019}; namely in these papers the sine-Gordon equation formulated in a semi-strip and the elliptic Ernst equation formulated in a domain representing the exterior of a thin rotating disk were analyzed. It is interesting that for these PDEs and these particular domains, there exist several linearizable boundary conditions. In recent years, the Fokas method has been used by many investigators for the analytical as well the numerical investigation of both linear and nonlinear elliptic PDEs, see for example \cite{AN2014, CFF2018, C2015, C2015b, DO2011, DT2014, FF2011, LC2018, OVDH2012, S2014, VD2013}.

\begin{figure}
       \vspace{.3cm}
\begin{center}
 \begin{overpic}[width=.5\textwidth]{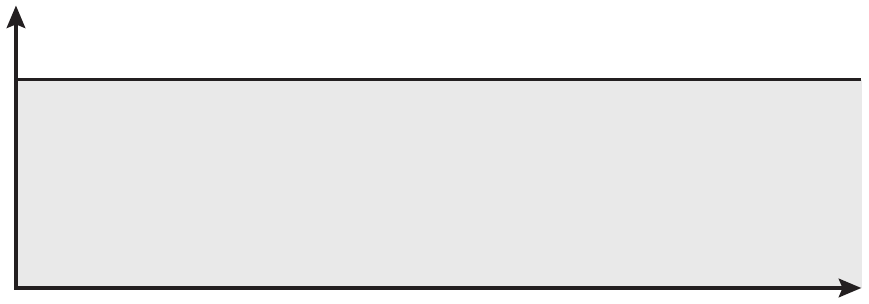}
      \put(101,.5){\small $x$}
      \put(0.5,37){\small $y$}
      \put(49.5,13){\small $\mathcal{S}$}
      \put(49,28){\tiny $(1)$}
      \put(-4.3,13.5){\tiny $(2)$}
      \put(49,-2.3){\tiny $(3)$}
       \put(-3,24.3){\small $L$}
      \end{overpic}
       \vspace{.3cm}
     \begin{figuretext}\label{semistrip.pdf}
       The semistrip $\mathcal{S}$ used in the formulation of the BVP for the elliptic sine-Gordon equation.
         \end{figuretext}
     \end{center}
\end{figure}

In this paper the following linearizable BVPs will be analyzed:

(a) Let $q(x,y)$ satisfy the elliptic sine-Gordon equation
\begin{equation} \label{esg}
  q_{xx}+q_{yy}=\sin q, 
\end{equation}
in the semistrip 
\begin{equation}
 \mathcal{S} = \{ 0<x<\infty, \quad 0<y<L\},
\end{equation}
depicted in Figure \ref{semistrip.pdf}, with the Dirichlet boundary conditions 
\begin{equation}\label{lbc}
q(x,0)=q(x,L)=0, \quad 0<x<\infty; \qquad q(0,y)=d, \quad 0<y<L, 
\end{equation}
where $L > 0$ and $d \in \R$ are finite constants.

\begin{figure}
       \vspace{.3cm}
\begin{center}
 \begin{overpic}[width=.4\textwidth]{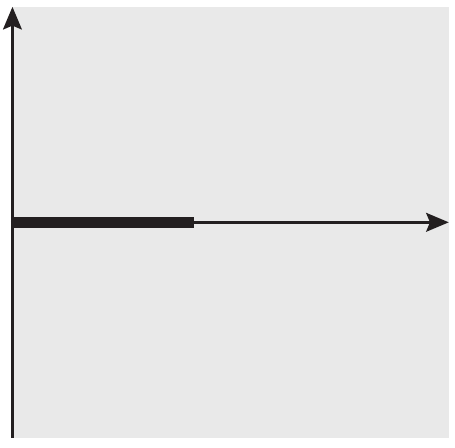}
      \put(0.5,100){\small $i\zeta$}
      \put(101,48.2){\small $\rho$}
      \put(40,44){\small $\rho_0$}
      \put(48,71){\small $\mathcal{D}$}
        \end{overpic}
     \begin{figuretext}\label{diskdomain.pdf}
       The exterior domain $\mathcal{D}$ of a finite disk of radius $\rho_0$ used in the formulation of the BVPs for the elliptic Ernst equation.
         \end{figuretext}
     \end{center}
\end{figure}

(b) Let $\mathcal{D} \subset \C$ denote the domain
$$\mathcal{D} = \{z \in \C \, | \, \re  z > 0\} \setminus [0, \rho_0],$$
where $\rho_0 > 0$ is a constant, see Figure \ref{diskdomain.pdf}. Let $z = \rho + i\zeta$ and consider the problem of finding a function $f(z)$ which satisfies the elliptic Ernst equation
\begin{equation}\label{ernst}  
  \frac{f + \bar{f}}{2}\left(f_{\rho\rho} + f_{\zeta\zeta} + \frac{1}{\rho} f_\rho\right) = f_\rho^2 + f_\zeta^2 \quad \text{in}\quad \mathcal{D},
\end{equation}  
together with the following boundary conditions:
\begin{subequations}\label{diskconditions}
\begin{align}\label{asymptoticflatness}
\bullet\ & \text{$f(z) \to 1$ as $|z|^2 \to \infty$ (asymptotic flatness),}
	\\ \label{axiscondition}
 \bullet \ & \text{$\frac{\partial f}{\partial \rho}(i\zeta) = 0$ for $|\zeta| > 0$  (regularity on the rotation axis),}
	\\ \label{BVPdisk}
\bullet \ &  \text{$f_{\Omega}(\rho \pm i 0) = \re  f_{\Omega}(+ i 0)$ for $0 < \rho < \rho_0$ (Dirichlet condition on the disk),}
\end{align}
\end{subequations}
where $\Omega \in \R$ is a constant and $f_\Omega$ denotes the Ernst potential in a frame corotating with angular velocity $\Omega$, see section \ref{ernstsec}.

(c) Let $f(z)$ satisfy the same BVP as in (b), but with the boundary condition (\ref{axiscondition}) replaced with
\begin{subequations}\label{holediskconditions}
\begin{align}
 \bullet \ & \text{$\frac{\partial f}{\partial \rho}(i\zeta) = 0$ for $|\zeta| > r_1$  (regularity on the rotation axis),}
	\\ \label{BVPhorizon}
 \bullet \ & \text{$\re  f_{\Omega_h}(i\zeta) = 0$ for $0 < |\zeta| < r_1$ (boundary condition on the horizon),}
\end{align}
\end{subequations}
where $\Omega_h \in \R$ and $r_1 > 0$ are constants.

(d) Let $f(z)$ satisfy the same BVP as in (b), but with the boundary condition (\ref{BVPdisk}) replaced with
\begin{align}
\bullet \ & \text{$\frac{\partial f_{\Omega}}{\partial \zeta}(\rho \pm i 0) = 0$ for $0 < \rho < \rho_0$ (Neumann condition on the disk).}
\end{align}

The BVPs (b) and (c) correspond to a uniformly rotating dust disk and a dust disk rotating uniformly around a central black hole, respectively, cf. \cite{MAKNP, L}. 
In fact, the solution of the BVP formulated in (b) is the celebrated Neugebauer-Meinel disk \cite{NM1993} of radius $\rho_0$ rotating with angular velocity $\Omega$. 
Moreover, if one sets $\rho_0 = 0$ in (c) (i.e. one removes the disk), then the solution of the obtained BVP is the Kerr black hole rotating with angular velocity $\Omega_h$. 
In the $z$-plane, the event horizon of this black hole stretches along the imaginary axis from $-ir_1$ to $ir_1$. 

The problems (a)-(d) were studied in detail using the unified transform in \cite{FLP2013}, \cite{LFernst}, \cite{L}, and \cite{LP2019}, respectively. The purpose of this paper is to review the approaches and results of these references while emphasizing the unified nature of the approach. 

\section{The sine-Gordon equation on the semistrip} \nequation
We will express the solution $q(x,y)$ in terms of the solution of a $2 \times 2$ RH problem (details are given in \cite{FLP2013}). As it was mentioned in the introduction, this RH problem depends on certain spectral functions.

\subsection{The spectral functions}
The spectral functions associated with the sides (1), (2), (3) indicated in Figure \ref{semistrip.pdf}, will be denoted respectively by
$$\{a_1(k), b_1(k)\}, \qquad \{a_2(k), b_2(k)\}, \qquad \{a_3(k), b_3(k)\}.$$
These functions are uniquely defined in terms of the following boundary values:
\begin{align*}
& \left\{g_0^{(1)}(x) = q(x, L), \;g_1^{(1)}(x) = q_y(x,L)\right\}, \qquad
\left\{g_0^{(2)}(y) = q(0, y), \; g_1^{(2)}(y) = q_x(0,y)\right\}, 
	\\
& \left\{g_0^{(3)}(x) = q(x, 0), \;g_1^{(3)}(x) = q_y(x,0)\right\}.
\end{align*}
We assume that
\begin{align}\nonumber
& g_0^{(1)}(x), \, g_1^{(1)}(x), \, g_0^{(3)}(x), \, g_1^{(3)}(x) \in L^1(\R^+),
	\\ \label{bcass}
& x g_0^{(1)}(x), \, x g_1^{(1)}(x), \, xg_0^{(3)}(x), \, x g_1^{(3)}(x) \in L^1(\R^+), 
	\\ \nonumber
& g_0^{(2)}(y), \, g_1^{(2)}(y), \, yg_0^{(2)}(y), \, yg_1^{(2)}(y) \in L^1([0, L]).
\end{align}
For a $2 \times 2$-matrix $A$, let $[A]_1$ and $[A]_2$ denote the first and second columns of $A$, respectively.
The functions $a_1(k)$ and $b_1(k)$ are defined by
\begin{subequations}\label{a1b1def}
\begin{equation} 
\begin{pmatrix} a_1(k) \\ b_1(k) \end{pmatrix}
= [\Phi_1(0, k)]_1, \qquad \im k \geq 0,
\end{equation}
where $\Phi_1(x,k)$ denotes the unique solution of the Volterra linear integral equation
\begin{equation}
  \Phi_1(x,k) = I - \int_x^\infty e^{\Omega(k) (\xi - x) \frac{\hat{\sigma}_3}{2}} Q^{(1)}(\xi, k) \Phi_1(\xi, k) d\xi,
\end{equation}
\end{subequations}
with $\Omega(k)$, $\hat{\sigma}_3$, and $Q^{(1)}$ defined as follows:
\begin{align} \nonumber
& \Omega(k) = \frac{1}{2i}\left(k - \frac{1}{k}\right);
	\\ \nonumber
& \hat{\sigma}_3A = [\sigma_3, A], \quad \sigma_3 = \text{diag}(1, -1), \quad \text{hence} \quad
e^{\hat{\sigma}_3}A = e^{\sigma_3}Ae^{-\sigma_3};
	\\ \nonumber
& Q^{(1)}(x, k) = \frac{i}{4}
\begin{pmatrix} \frac{1}{k}(1 - \cos g_0^{(1)}) & \dot{g}_0^{(1)} - i g_1^{(1)} + \frac{i}{k} \sin g_0^{(1)} \\
\dot{g}_0^{(1)} - i g_1^{(1)} - \frac{i}{k} \sin g_0^{(1)} & - \frac{1}{k}(1 - \cos g_0^{(1)}) 
\end{pmatrix}, 
	\\ \nonumber
& \dot{g}_0^{(1)} = \frac{d}{dx}g_0^{(1)}(x).
\end{align}	
The functions $a_2(k)$ and $b_2(k)$ are defined by
\begin{subequations}\label{a2b2def}
\begin{equation} 
\begin{pmatrix} a_2(k) \\ b_2(k) \end{pmatrix}
= [\Phi_2(0, k)]_1, \qquad \im k \geq 0,
\end{equation}
where $\Phi_2(y,k)$ denotes the unique solution of the Volterra linear integral equation
\begin{equation}
  \Phi_2(y,k) = I - i \int_y^L e^{\omega(k) (\eta - y) \frac{\hat{\sigma}_3}{2}} Q^{(2)}(\eta, k) \Phi_2(\eta, k) d\eta,
\end{equation}
\end{subequations}
with $\omega(k)$ and $Q^{(2)}$ defined as follows:
\begin{align} \nonumber
& \omega(k) = \frac{1}{2}\left(k + \frac{1}{k}\right);
	\\ \nonumber
& Q^{(2)}(y, k) = \frac{i}{4}\begin{pmatrix} \frac{1}{k}(1 - \cos g_0^{(2)}) & g_1^{(2)} - i \dot{g}_0^{(2)} + \frac{i}{k} \sin g_0^{(2)} \\ \nonumber
g_1^{(2)} - i \dot{g}_0^{(2)} - \frac{i}{k} \sin g_0^{(2)} & - \frac{1}{k}(1 - \cos g_0^{(2)}) 
\end{pmatrix}, 
	\\ \nonumber
& \dot{g}_0^{(3)} = \frac{d}{dy}g_0^{(3)}(y).
\end{align}	
The functions $a_3(k)$ and $b_3(k)$ are defined by equations which are similar to (\ref{a1b1def}) but with $Q^{(1)}$ replaced by $Q^{(3)}$.

\subsection{The global relation}
The spectral functions satisfy the following pair of global relations:
\begin{subequations}\label{globalrelations}
\begin{align}\label{globalrelation1}
& a_1(k) = a_2(-k)a_3(k) - b_2(-k)b_3(k),  & k \in \C^+,
	\\ \label{globalrelation2}
& b_1(k)e^{-\omega(k)L}  = a_2(k)b_3(k) - a_3(k)b_2(k),  & k \in \C^+.
\end{align}
\end{subequations}

\begin{theorem}{\upshape \bf \cite{FLP2013}} \label{mainRHth}
Suppose that a subset of the boundary values $\{q(x,L), q_y(x,L)\}$, $\{q(x,0), q_y(x,0)\}$, $0<x<\infty$, and $\{q(y,0), q_x(y,0)\}$, $0<y<L$, satisfying (\ref{bcass}), are prescribed as boundary conditions. Suppose that these prescribed boundary conditions are such that the global relations (\ref{globalrelation1}) and (\ref{globalrelation2}) can be used to characterize the remaining boundary values.
Define the spectral functions $\{a_j, b_j\}$, $j=1, 2,3,$ by (\ref{a1b1def})-(\ref{a2b2def}). Assume that $a_1(k)$ and $a_3(k)$ do not have zeros for $\text{\upshape Im}\, k \geq 0$.

Define $M(x,y,k)$ as the solution of the following $2\times 2$ matrix RH problem:

\begin{itemize}
\item The function $M(x,y,k)$ is a sectionally meromorphic function of $k$ away from $\R\cup i\R$.

\item $M=I+O\left(\frac 1 k\right)$ as $k\to\infty$.

\item $M$ satisfies the jump condition
$$M_-(x,y,k)=M_+(x,y,k)J(x,y,k),\qquad k\in\R\cup i\R,$$
where $M=M_+$ for $k$ in the first or third quadrant, and $M=M_-$ for $k$ in the second or fourth quadrant of the complex $k$ plane, and $J$ is defined in terms of $\{a_j,b_j\}$ as follows:
\begin{align}\nonumber
& J(x,y,k) =
J^{\alpha}(x,y,k),\qquad \text{\upshape arg}\,k=\alpha,\quad
\alpha=0,\,\frac{\pi}{2},\, \pi,\,\frac{3\pi}{2};
	\\ \label{alljumps2}
& J^{\pi/2}=\begin{pmatrix}
1 & I(k)e^{-\theta(x,y,k)} 	\\ 
0 & 1
\end{pmatrix}, \qquad
J^{3\pi/2}=\begin{pmatrix}
1& 0 \\
I(-k)e^{\theta(x,y,k)}&1
 \end{pmatrix},
	\\  \nonumber
& J^{0} = \begin{pmatrix}
R(k)& \frac {b_3(-k)}{a_3(k)}e^{-\theta(x,y,k)}
\\ -\frac{e^{-\omega(k)L}b_1(k)}{a_1(-k)}e^{\theta(x,y,k)}& 1 \end{pmatrix},
\qquad J^{\pi} = J^{3\pi/2}(J^{0})^{-1}J^{\pi/2},
\end{align}
where the functions $\theta, I, R$ are defined by
\begin{align}\nonumber
& \theta(x,y,k)=\Omega(k)x+\omega(k)y, \qquad
 I(k) = \frac{b_2(-k)}{a_1(k)a_3(k)}, 
	\\\label{Theta}
& R(k) = \frac {a_2(k)}{a_1(-k)a_3(k)}.
\end{align}
\end{itemize}

Then, $M$ exists and is unique, provided that the $H^1$ norm of the spectral functions is sufficiently small.

Define $q(x,y)$ in terms of $M(x,y,k)$ by
\begin{subequations}\label{qrecover}
\begin{align}
q_{x}-iq_y = &\; 2\lim_{k\to\infty}(k M)_{21},
 	\\
  \cos q(x,y) = & \; 1+4i\Bigl(\lim_{k\to\infty}(k M_x)_{11}\Bigr) - 2\Bigl(\lim_{k\to\infty}(k M)_{21}\Bigr)^2.
\end{align}                         
\end{subequations}           
Then $q(x,y)$ solves (\ref{esg}). Furthermore, 
\begin{align*}
 & q(x, L) = g_0^{(1)}(x), \quad q_y(x,L) = g_1^{(1)}(x), 
 	\\
& q(x, 0) = g_0^{(3)}(x), \quad q_y(x,0) = g_1^{(3)}(x), \quad 0 <  x < \infty;
	\\
&q(0, y) = g_0^{(2)}(y), \quad  q_x(0,y) = g_1^{(2)}(y), \qquad 0 < y < L.
\end{align*}
\end{theorem}

\section{Linearizable boundary conditions for the elliptic sine-Gordon equation}\nequation
We now concentrate on the particular boundary conditions (\ref{lbc}).
In this case, using the notations
\begin{align*}
& \Phi_1(x,k) = \begin{pmatrix} A_1(x,k) & B_1(x, -k) \\  B_1(x,k) & A_1(x,-k) \end{pmatrix}, \qquad
\Phi_2(y,k) = \begin{pmatrix} A_2(y,k) & B_2(y, -k) \\  B_2(y,k) & A_2(y,-k) \end{pmatrix},
	\\
& \Phi_3(x,k) = \begin{pmatrix} A_3(x,k) & B_3(x, -k) \\  B_3(x,k) & A_3(x,-k) \end{pmatrix},
\end{align*}
the integral equations characterising the spectral functions simplify as follows:
\begin{subequations}
\begin{align} \label{A1eq}
& \begin{pmatrix}
A_1(x,k)\\B_1(x,k)
\end{pmatrix}
=\begin{pmatrix} 1 \\ 0 \end{pmatrix}
- \frac{1}{4} \int_x^{\infty}\begin{pmatrix}
q_y(\xi,L)B_1(\xi,k)\\ e^{\Omega(k)(x-\xi)}q_y(\xi,L)A_1(\xi,k) \end{pmatrix}d\xi, \quad \; 0<x<\infty,\ \im k \geq 0,
	\\ \nonumber 
&\begin{pmatrix} A_2(y,k)\\B_2(y,k)\end{pmatrix}
=\begin{pmatrix} 1 \\ 0 \end{pmatrix}
+\frac{1}{4}\int_y^{L}\begin{pmatrix}
-\frac{(1-\cos d)}{k}A_2(\eta,k)+[q_x(0,\eta)-i\frac{\sin d}{k}]B_2(\eta,k) \\
e^{\omega(k)(y-\eta)}[q_x(0, \eta)+i\frac{\sin d}{k}]A_2(\eta,k)+
\frac{(1-\cos d)}{k}B_2(\eta,k)\end{pmatrix} d\eta,
	\\ \label{A2eq}
&\hspace{9cm}  0<y<L,\ k\in\C,
	\\  \label{A3eq}
& \begin{pmatrix}
A_3(x,k)\\B_3(x,k)\end{pmatrix} 
=\begin{pmatrix} 1\\0\end{pmatrix}
-\frac 1 4\int_x^{\infty}\begin{pmatrix}
q_y(\xi,0)B_3(\xi,k)\\ e^{\Omega(k)(x-\xi)}q_y(\xi,0)A_3(\xi,k)
\end{pmatrix} d\xi, \quad \;  0<x<\infty,\   \im k \geq 0.
\end{align}
\end{subequations}
In equations (\ref{A1eq}) and (\ref{A3eq}), the only dependence on $k$ is through $\Omega(k)$.  Thus, since $\Omega(-1/k)=\Omega(k)$, it follows that the vector functions $(A_1,B_1)$ and $(A_3,B_3)$ satisfy the same symmetry properties. Hence,
\begin{equation}\label{symm1}
  a_j\Bigl(-\frac 1 k\Bigr)=a_j(k),\quad b_j \Bigl(-\frac 1 k \Bigr)=b_j(k),\qquad j=1,3, \ \im k\geq 0.
\end{equation}
It turns out that the vector function $(A_2,B_2)$ also satisfies a certain symmetry condition, as stated in the following proposition.
\begin{proposition}
Let $q_x(0,y)$ be a sufficiently smooth function. Then the vector solution of the linear Volterra integral equation (\ref{A2eq}) satisfies the following symmetry conditions for $0<y<L$ and $k\in\C$:
\begin{align} \label{symm2}
\begin{cases}
 A_2\big(\frac 1 k\big)=\frac{A_2(k)-F(k)B_2(k)+F(k)e^{\omega(k)(y-L)}B_2(-k)-
F(k)^2e^{\omega(k)(y-L)}A_2(-k)}{1-F(k)^2},
\vspace{.1cm}	 \\ 
 B_2\big(\frac 1 k \big)=\frac{B_2(k)-F(k)A_2(k)+F(k)e^{\omega(k)(y-L)}A_2(-k)-
F(k)^2e^{\omega(k)(y-L)}B_2(-k)}{1-F(k)^2},
\end{cases}
\end{align}
where the $y$-dependence of $A_2$ and $B_2$ has been suppressed for clarity and 
\begin{equation}\label{Fdef}
  F(k)=i\frac{1-k^2}{1+k^2}\tan\frac d 2.
\end{equation}
\end{proposition}

Recalling that $a_2(k)=A_2(0,k)$, and $b_2(k)=B_2(0,k)$, equations (\ref{symm2}) immediately imply the following
important relations:
\begin{align} \label{symm3}
\begin{cases}
a_2\big(\frac 1 k\big)=\frac{a_2(k)-F(k)b_2(k)+F(k)e^{-\omega(k)L}b_2(-k)-
F(k)^2e^{-\omega(k)L}a_2(-k)}{1-F(k)^2},
 \vspace{.1cm}	\\ 
b_2 \big(\frac 1 k \big) =\frac{b_2(k)-F(k)a_2(k)+F(k)e^{-\omega(k)L}a_2(-k)-
F(k)^2e^{-\omega(k)L}b_2(-k)}{1-F(k)^2},
\end{cases} \quad \im k \geq 0.
\end{align}
In summary, the basic equations characterizing the spectral functions are:
\begin{itemize}
\item[(a)] the symmetry relations (\ref{symm1}) and (\ref{symm3});
\item[(b)] the global relations (\ref{globalrelation1}) and (\ref{globalrelation2});
\item[(c)] the conditions of unit determinant.
\end{itemize}
It turns out that, using these equations, it is possible to provide an explicit characterization of all the spectral functions in terms of the given constant $d$.

\begin{proposition}
Assume that the functions $\{a_j(k),b_j(k)\}$, $j=1,2,3$ satisfy the symmetry relations (\ref{symm1}) and (\ref{symm3}), the global relations (\ref{globalrelation1}) and (\ref{globalrelation2}), and the unit determinant conditions
\begin{equation}\label{dc}
  a_j(k)a_j(-k)-b_j(k)b_j(-k)=1,\qquad j=1,2,3.
\end{equation}
Then the following relations are valid:
\begin{subequations}
\begin{align}\label{rel1}
&a_1(k)b_1(-k)-a_1(-k)b_1(k)=G(k),&\qquad k\in\R,
  	\\ \label{rel2}
&a_3(k)b_3(-k)-a_3(-k)b_3(k)=-G(k),&\qquad k\in\R,
	\\ \label{rel3}
&a_2(k)=a_1(-k)a_3(k)-e^{-\omega(k)L}b_1(k)b_3(-k),&\qquad k\in\C,
	\\ \label{rel4}
&b_2(k)=a_1(-k)b_3(k)-e^{-\omega(k)L}b_1(k)a_3(-k),&\qquad k\in\C,
\end{align}
\end{subequations}
where
\begin{equation}\label{Gdef}
G(k)=\frac{i(1-k^2)}{1+k^2}\frac {e^{\omega(k)L}+e^{-\omega(k)L}-2}{e^{\omega(k)L}-e^{-\omega(k)L}}
\tan\frac d 2.
\end{equation}
\end{proposition}

\begin{remark} The two relations (\ref{rel3}) and (\ref{rel4}) are a direct consequence of the global relation and of the conditions of unit determinant. On the other hand, equations (\ref{rel1}) and (\ref{rel2}) depend  on the particular symmetry properties.
\end{remark}

\begin{remark}
The determinant condition (\ref{dc}) with $j=1$ and equation (\ref{rel1}) imply
\begin{equation}\label{fact0}
\left[a_1(k)^2-b_1(k)^2\right]\left[a_1(-k)^2-b_1(-k)^2\right]=1-G(k)^2,\qquad k\in\R.
\end{equation}
Indeed, equation (\ref{dc}) with $j=1$ implies the identity
\begin{equation}\label{fact01}
\bigl(a_1^2-b_1^2\bigr)\bigl(\hat a_1^2-\hat b_1^2\bigr) = 1- \bigl(a_1\hat b_1-\hat a_1b_1\bigr)^2,\qquad k\in\R,
\end{equation}
where, for a scalar-valued function $f$, we use the notation $\hat{f} := f(-k)$.
Then equation (\ref{rel1}) implies (\ref{fact0}).
\end{remark}

Equation (\ref{fact0}) defines the jump relation of a scalar RH problem for the sectionally analytic function defined by
$$\begin{cases}
a_1(k)^2-b_1(k)^2, & \qquad k\in\C^+, \\ 
 a_1(-k)^2-b_1(-k)^2,& \qquad k\in\C^-.
\end{cases}
$$
Taking into consideration that $a_1(k)\neq b_1(k)$ for $k\in\C^+$ (otherwise equation (\ref{dc}) with $j=1$ is violated) it follows that the above Riemann--Hilbert problem has a unique solution
\begin{equation}\label{sol1}
a_1(k)^2-b_1(k)^2=h(k),\qquad k\in\C^+,
\end{equation}
where
\begin{equation}\label{rhH}
h(k)= e^{H(k)}, \quad H(k)=\frac 1 {2\pi i}\int_{\R}\ln[1-G^2(k')]\frac{dk'}{k'-k},\qquad k\in\C\setminus\R.
\end{equation}
Using the fact that $G(k)$ is an odd function, it follows that $H(k)$ is also an odd function, hence
$h(-k) = e^{-H(k)}$. This implies that the function $h(k)$ defined by (\ref{rhH}) satisfies the jump condition (\ref{fact0}).

\subsection{Spectral theory in the linearizable case}
In the case of  the linearizable boundary conditions (\ref{lbc}), it is possible to express $q(x,y)$ in terms of the solution of a RH problem whose jump matrices are computed explicitly in terms of the given constant $d \in \R$. Indeed, equations (\ref{rel3}) and (\ref{rel4}) imply that the functions $R(k)$ and $I(k)$ in (\ref{Theta}) are given by
\begin{equation}\label{IR}
R(k)=1-e^{-\omega(k)L}\frac{\hat b_3}{a_3}\frac {b_1}{\hat a_1},
\qquad
I(k)=\frac {\hat b_3}{a_3}-e^{\omega(k)L}\frac{\hat b_1}{a_1}.
\end{equation}
Thus in the linearizable case, the jump matrices involve only the ratios $\frac {\hat b_3}{a_3}$ and $\frac{\hat b_1}{a_1}$, evaluated at $k$ and at $-k$. Equations (\ref{sol1}) and (\ref{rhH}) imply that these ratios are given by
\begin{equation}\label{hatnohat}
\frac {\hat b_3}{a_3}=-\frac{G}{h}+\frac{b_3}{a_3h},\qquad
\frac {\hat b_1}{a_1}=\frac{G}{h}+\frac{b_1}{a_1h}.
\end{equation}
Hence the jump matrices depend on the known function $G/h$ as well as on the unknown functions
$\frac{b_1}{a_1h}$ and $\frac{b_3}{a_3h}$. Using the fact that these unknown functions are bounded and analytic in $\C^+$, it is possible to formulate a RH problem in terms of the known function $G/h$ alone, which is equivalent to the basic RH problem defined in Theorem \ref{mainRHth}.

\begin{theorem}{\upshape \bf \cite{FLP2013}}
Let $q(x,y)$ satisfy equation (\ref{esg}) and the boundary conditions (\ref{lbc}).
Then $q(x,y)$ is given by equations (\ref{qrecover}) with $M$ replaced by $\tilde{M}$, where $\tilde{M}$ is the solution of the Riemann--Hilbert problem of Theorem \ref{mainRHth} with the jump matrix $J$ replaced by the matrix $\tilde J$  defined as follows:
\begin{align*}
& \tilde J(x,y,k)=
\tilde J^{\alpha}(x,y,k), \qquad \text{\upshape arg}\, k = \alpha,\qquad
\alpha=0,\,\frac{\pi}{2},\, \pi,\,\frac{3\pi}{2};
	\\
& \tilde J^{0}=\begin{pmatrix}
1-\frac{G^2(k)}{h(k) h(-k)}e^{-\omega(k)L}&-\frac{G(k)}{h(k)}e^{-\theta(x,y,k)}
\\
e^{-\omega(k)L}\frac{G(k)}{h(-k)}e^{\theta(x,y,k)}& 1
 \end{pmatrix},
	\\ \label{alljumpslin}
& \tilde J^{\pi/2}=\begin{pmatrix}
1&
-\frac{G(k)}{h(k)}\left(1+e^{\omega(k)L}\right)e^{-\theta(x,y,k)}
\\
0&1
 \end{pmatrix},
 	\\
& \tilde  J^{3\pi/2}=\begin{pmatrix}
1& 0
\\
\frac{G(k)}{h(-k)}\left(1+e^{-\omega(k)L}\right)e^{\theta(x,y,k)}&1
 \end{pmatrix}, \qquad
\tilde J^{\pi}=
\tilde J^{3\pi/2}(\tilde J^{0})^{-1}\tilde J^{\pi/2},
\end{align*}
where $\theta(x,y,k) = \Omega(k)x+\omega(k)y$, while $G(k)$ and $h(k)$ are defined in terms of the given constant $d$ by equations (\ref{Gdef}) and (\ref{rhH}).

This RH problem is regular and has a unique solution.
\end{theorem}

\section{The elliptic Ernst equation in an exterior disk domain}\label{ernstsec}\nequation
We first consider the case of a rotating disk in the absence of a central black hole. Thus, let $f(z)$ be a solution of the elliptic Ernst equation (\ref{ernst}) in $\mathcal{D}$, which satisfies the boundary conditions (\ref{asymptoticflatness}) and (\ref{axiscondition}) together with some compatible, but otherwise arbitrary, boundary conditions along the disk.

For each $z \in \mathcal{D}$, we let $\mathcal{S}_z$ denote the compact Riemann surface of genus $0$ defined by the equation
\begin{equation*}
  \lambda = \sqrt{\frac{k - i \bar{z}}{k + i z}}.
\end{equation*}
We view $\mathcal{S}_z$ as a two-sheeted covering of the Riemann $k$-sphere endowed with a branch cut from $-iz$ to $i\bar{z}$; the upper (lower) sheet is characterized by $\lambda \to 1$ ($\lambda \to -1$) as $k \to \infty$. We let $k^+$ and $k^-$ denote the points which project onto $k \in \hat{\C} = \C \cup \{\infty\}$ and which lie in the upper and lower sheet of $\mathcal{S}_z$, respectively.
The next proposition expresses $f(z)$ in terms of the solution of a $2\times2$ matrix RH problem on $\mathcal{S}_z$ formulated in terms of {\it both} the Dirichlet and the Neumann boundary values on the disk.

\begin{proposition}\label{RHprop}
 Let $f(z)$ satisfy the Ernst equation (\ref{ernst}) in the exterior disk domain $\mathcal{D}$. Suppose that $f$ satisfies the boundary conditions (\ref{asymptoticflatness}) and (\ref{axiscondition}) and that $\text{\upshape Re} f > 0$ in $\mathcal{D}$. The solution $f(z)$ can be expressed in terms of the Dirichlet and the Neumann boundary values on the disk as follows:
\begin{enumerate}
\item Use the disk boundary values of $f$ to define the $2\times 2$-matrix valued one-forms $W^{R}(\rho \pm i0, k)$ and $W^{L}(\rho \pm i0, k)$ for $0 < \rho < \rho_0$ and $k = ik_2$, $0 < k_2 < \rho_0$, by
\begin{equation}\label{Wrightorabovecut}
W^{R}(\rho \pm i0, k)
= \begin{cases}
\frac{1}{f + \bar{f}} \begin{pmatrix} \bar{f}_\rho & \frac{ik_2 \bar{f}_\rho - \rho \bar{f}_\zeta}{i\sqrt{k_2^2 - \rho^2}} \\
\frac{i k_2 f_\rho - \rho f_\zeta}{i\sqrt{k_2^2 - \rho^2}} & f_\rho 
\end{pmatrix}_{z = \rho \pm i0} d\rho, \qquad 0 < \rho < k_2,
	\\
\frac{1}{f + \bar{f}} \begin{pmatrix} \bar{f}_\rho & \frac{ik_2 \bar{f}_\rho - \rho \bar{f}_\zeta}{\sqrt{\rho^2 - k_2^2}} \\
\frac{ik_2 f_\rho - \rho f_\zeta}{\sqrt{\rho^2 - k_2^2}} & f_\rho 
\end{pmatrix}_{z = \rho \pm i0} d\rho, \qquad k_2 < \rho < \rho_0
\end{cases}
\end{equation}
and
\begin{equation}\label{WLR}
W^{L}(\rho \pm i0, k)
= \begin{cases}
W^{R}(\rho \pm i0, k), \qquad 0 < \rho < k_2, \\
\sigma_3 W^{R}(\rho \pm i0, k) \sigma_3, \qquad k_2 < \rho < \rho_0.
\end{cases}
\end{equation}

\item Define a function $\Psi(\rho \pm i0, k)$ for $0 < \rho<\rho_0$ and $k = ik_2$, $0 < k_2 < \rho_0$, by solving the ordinary differential equation
$$\Psi_\rho(\rho \pm i0, k) = W^L(\rho \pm i0, k)\Psi(\rho \pm i0, k), \qquad 0 < \rho < \rho_0,$$
together with the initial conditions
$$[\Psi(-i0, k)]_1 = \begin{pmatrix} \overline{f(-i0)} \\ f(-i0) \end{pmatrix}, \qquad
[\Psi(+i0, k)]_2 = \begin{pmatrix} 1 \\ -1 \end{pmatrix},$$
as well as the following continuity condition at the rim of the disk:
$$\Psi(\rho_0 - i0, k) = \Psi(\rho_0 + i0, k).$$

\item Let
$$C(k) = \int_{[k_2, \rho_0]} \left[\left((W^R - W^L)\Psi\right)(\rho - i0, k) - \left((W^R - W^L)\Psi\right)(\rho + i0, k) \right],$$
and define the functions $F^\pm(k), G^\pm(k)$ by solving the algebraic system
$$\begin{pmatrix} 1 & G_-(k) \\ 0 & F_-(k) \end{pmatrix} - \begin{pmatrix} 1 & G_+(k) \\ 0 & F_+(k) \end{pmatrix}\begin{pmatrix} F_+(k) & 0 \\ G_+(k) & 1 \end{pmatrix}^{-1}\begin{pmatrix} F_-(k) & 0 \\ G_-(k) & 1 \end{pmatrix}
= \frac{1}{2}\begin{pmatrix} 1 & 1 \\ 1 & - 1 \end{pmatrix}C(k).$$

\item Define a $2\times 2$-matrix valued function $D(k)$, $k = ik_2$, $0 < k_2 < \rho_0$ by
\begin{equation}\label{Dequation1}
  D(k) = \begin{pmatrix} F_+(k) & 0 \\ G_+(k) & 1 \end{pmatrix}^{-1}\begin{pmatrix} F_-(k) & 0 \\ G_-(k) & 1 \end{pmatrix}.
\end{equation} 
The entries of $D(k)$ are rational functions of the entries of $C(k)$. Define $D(ik_2)$ for $-\rho_0 < k_2 < 0$ by $D(ik_2) = \sigma_3 \overline{D(-ik_2)} \sigma_3$.

\item Let $\Phi(z, k)$ be the unique solution of the following RH problem:
\begin{itemize}
\item For each $z$, $\Phi(z, \cdot)$ is a map from the Riemann surface $\mathcal{S}_z$ to the space of $2 \times 2$ matrices.

\item $\Phi(z, k)$ is an analytic function of $k \in \mathcal{S}_z \setminus (\Gamma^+ \cup \Gamma^-)$, where $\Gamma^+$ and $\Gamma^-$ denote the coverings of 
$\Gamma = [-i\rho_0, i\rho_0]$ in the upper and lower sheets of $\mathcal{S}_z$, respectively.

\item Across $\Gamma^+$, $\Phi$ satisfies the jump condition 
$$ \Phi_-(z, k) = \Phi_+(z, k) D(k), \qquad k \in \Gamma^+,$$
where $\Phi_+$ and $\Phi_-$ denote the values of $\Phi$ infinitesimally to the right and left of $\Gamma^+$, respectively.

\item Across $\Gamma^-$, $\Phi$ satisfies the jump condition 
$$ \Phi_-(z, k) = \Phi_+(z, k) \sigma_1 D(k) \sigma_1, \qquad k \in \Gamma^-, \quad \text{where} \quad \sigma_1 = \begin{pmatrix} 0	&	1 \\ 1	& 0\end{pmatrix}.$$

\item As $k \to \infty$, $\Phi$ satisfies
\begin{equation}\label{RHnormalization} 
\lim_{k \to \infty} [\Phi(z, k^-)]_1 = \begin{pmatrix} 1 \\ 1\end{pmatrix}, \qquad \lim_{k \to \infty} [\Phi(z, k^+)]_2 = \begin{pmatrix} 1 \\ -1\end{pmatrix}.
\end{equation}

\item $\Phi$ obeys the symmetries
\begin{equation}\label{PhisymmRH} 
\Phi(z, k^+) = \sigma_3 \Phi(z, k^-) \sigma_1, \qquad 
\Phi(z, k^+) = \sigma_1 \overline{\Phi(z, \bar{k}^+)} \sigma_3.
\end{equation}
\end{itemize}

\item Find $f(z)$ from the equation
\begin{equation}\label{frecover}
  f(z) = \lim_{k \to \infty} (\Phi(z, k^+))_{21}.
\end{equation}  
\end{enumerate}
\end{proposition}

\subsection{Equatorial symmetry and the global relation}
Since only a subset of the boundary values can be specified for a well-posed problem, the solution formula presented in Proposition \ref{RHprop} is {\it not} yet effective. However, for equatorially symmetric solutions whose boundary values possess a sufficient amount of symmetry (such boundary values are called linearizable), the unknown boundary values can be eliminated by using the global relation. 

The elliptic Ernst equation (\ref{ernst}) admits the Lax pair
\begin{equation}
\begin{cases}\label{ernstlax}
\Phi_z(z, k) = U(z, k) \Phi(z, k),
	\\
\Phi_{\bar{z}}(z, k) =V(z,k) \Phi(z,k),
\end{cases}
\end{equation}
where the $2 \times 2$-matrix valued function $\Phi(z, k)$ is an eigenfunction, $k$ is a spectral parameter, and the $2\times 2$-matrix valued functions $U$ and $V$ are defined by
$$U = \frac{1}{f + \bar{f}}\begin{pmatrix} \bar{f}_z & \lambda \bar{f}_z \\
\lambda f_z & f_z\end{pmatrix}, \qquad 
V =  \frac{1}{f + \bar{f}}\begin{pmatrix} \bar{f}_{\bar{z}} & \frac{1}{\lambda} \bar{f}_{\bar{z}} \\
\frac{1}{\lambda}  f_{\bar{z}}   & f_{\bar{z}}  \end{pmatrix}, 
\qquad \lambda(z,k) = \sqrt{\frac{k - i\bar{z}}{k + iz}}.$$
We let $\Phi$ be the solution of (\ref{ernstlax}) which satisfies the initial conditions
\begin{align}\label{phiinitial}
\lim_{z \to i\infty} [\Phi(z, k^-)]_1 = \begin{pmatrix} 1 \\ 1 \end{pmatrix}, \qquad
  \lim_{z \to i\infty} [\Phi(z, k^+)]_2 = \begin{pmatrix} 1 \\ -1 \end{pmatrix}, \qquad k \in \hat{\C} = \C \cup \infty.
\end{align}

\begin{proposition}\label{diskPhionaxisprop}
The values of $\Phi$ on the rotation axis $\rho = 0$ can be expressed in terms of two spectral functions $F(k)$ and $G(k)$, $k \in \hat{\C}$, as follows:
\begin{align*} 
\Phi(i\zeta, k^+) =  \begin{pmatrix} \overline{f(i\zeta)} & 1 \\ f(i\zeta) & -1 \end{pmatrix}\begin{pmatrix} F(k) & 0 \\ G(k) & 1 \end{pmatrix}, \qquad \zeta  > 0,\quad k \in \hat{\C},
	\\ 
\Phi(i\zeta, k^+) = \begin{pmatrix} \overline{f(i\zeta)} & 1 \\ f(i\zeta) & -1 \end{pmatrix}\begin{pmatrix} 1 & G(k) \\ 0 & F(k) \end{pmatrix}, \qquad \zeta  < 0,\quad k \in \hat{\C}.
\end{align*}
\end{proposition}

\begin{proposition}\label{eqprop}
Suppose that $f$ is equatorially symmetric, i.e., that $f(z) = \overline{f(\bar{z})}$ for $z \in \mathcal{D}$. 
Then the spectral functions $F(k)$ and $G(k)$ satisfy the following relation, which will be referred to as the global relation:
\begin{equation}\label{eqrelation}
  \overline{A_+(k)}\sigma_1\overline{A_+^{-1}(k)} \sigma_1 = \sigma_1A_-(k) \sigma_1A_-^{-1}(k), \qquad k \in \Gamma,
\end{equation}
where
\begin{equation}\label{Adef}
A(k) := \begin{pmatrix} F(k) & 0 \\ G(k) & 1 \end{pmatrix},
\end{equation}
and $A_\pm$ denote the values of $A$ to the right and left of $\Gamma$, respectively. 
\end{proposition}

\begin{remark}
For the BVPs denoted by {\upshape (b)-(d)} in the introduction, symmetry considerations and the assumption of uniqueness imply that the solution is equatorially symmetric.
\end{remark}

\section{Linearizable boundary conditions for the elliptic Ernst equation}\nequation
For a linearizable BVP, the spectral functions $F(k)$ and $G(k)$ satisfy, in addition to the global relation (\ref{eqrelation}), an additional important algebraic relation. 
These two algebraic relations satisfied by $F(k)$ and $G(k)$ yield an auxiliary RH problem for $F$ and $G$ with jump data formulated in terms of the known boundary values alone.

\subsection{A rotating disk}
Suppose $f(z)$ is a solution of the BVP formulated in (b) of the introduction. We will show that the boundary conditions satisfied by $f$ are linearizable and derive an explicit expression for the solution $f$ in terms of theta functions (see \cite{LFernst} for details). In this way we recover the Neugebauer-Meinel disk solutions \cite{NM1993}.

The boundary conditions (\ref{diskconditions}) and (\ref{holediskconditions}) involve the corotating Ernst potentials $f_\Omega$ and $f_{\Omega_h}$, which are defined as follows.
Outside a stationarily rotating axisymmetric body, the Einstein field equations are equivalent to the elliptic Ernst equation (\ref{ernst}). Indeed, in canonical Weyl coordinates the exterior gravitational field of such a body is described by the line element
\begin{equation}\label{lineelement}  
  ds^2 = e^{-2U}\bigl[e^{2\kappa}(d\rho^2 + d\zeta^2) + \rho^2d\varphi^2\bigr] - e^{2U}(dt + a d\varphi)^2,
\end{equation}
where $\rho, \zeta, \varphi$ are cylindrical coordinates, $t$ is the coordinate time, and the metric functions $e^{2U}, a, e^{2\kappa}$ depend only on $\rho$ and $\zeta$. 
Letting $f = e^{2U} + ib$, where $a_z = i \rho e^{-4U} b_z$, Einstein's equations reduce to (\ref{ernst}). 
Given $\Omega \in \R$, we define the coordinates $(\rho', \zeta', \varphi', t')$ corotating with the angular velocity $\Omega$ by
$$\rho' = \rho, \qquad \zeta' = \zeta, \qquad \varphi' = \varphi - \Omega t, \qquad t' = t.$$
The Ernst equation retains its form in the corotating system and $f_\Omega$ denotes the corresponding Ernst potential.

\begin{proposition}\label{fconstprop}
Suppose that $f_\Omega$ is constant along the disk. Let $f_0 := f(+i0)$ and define $B$ and $\Lambda(k)$ by 
\begin{equation}\label{Bdef}
  B = \begin{pmatrix} \overline{f_0} &  1  \\	 f_0 &   -1 \end{pmatrix}, \qquad \Lambda(k) = I + \frac{i k \Omega}{\text{\upshape Re}\, f_0}(\sigma_1 - I)\sigma_3.
\end{equation}
Then the spectral functions $F(k)$ and $G(k)$ satisfy the relation
\begin{equation*}
  (B^{-1}\Lambda^{-1}\sigma_1\sigma_3\overline{\Lambda B})(\overline{A_+ \sigma_1 A_+^{-1}}) = - (A_+\sigma_1A_+^{-1})(B^{-1} \Lambda^{-1}\sigma_1\sigma_3 \overline{\Lambda B}), \quad k \in \Gamma.
\end{equation*}
\end{proposition}

Combining Propositions \ref{eqprop} and \ref{fconstprop} we can determine the spectral functions $F$ and $G$ via the solution of a $2\times 2$ matrix RH problem.

\begin{proposition}\label{auxRHprop}
Suppose $f$ satisfies the BVP denoted by (b) in the introduction. 
Then the spectral functions $F(k)$ and $G(k)$ are given by
$$F(k) = -\mathcal{M}_{21}(k), \qquad G(k) = \mathcal{M}_{11}(k), \qquad k \in \C,$$
where $\mathcal{M}$ is the unique solution of the following RH problem:
\begin{itemize}
\item \text{$\mathcal{M}(k)$ is analytic for $k \in \C \backslash \Gamma$, $\Gamma = [-i\rho_0, i\rho_0]$.} 

\item  Across $\Gamma$, $\mathcal{M}(k)$ satisfies the jump condition
\begin{align}\label{calMjump}  
   \mathcal{S}(k)\mathcal{M}_-(k) = -\mathcal{M}_+(k)\mathcal{S}(k),\qquad k \in \Gamma, 
\end{align}
where $\mathcal{M}_+$ and $\mathcal{M}_-$ denote the values of $\mathcal{M}$ to the right and left of $\Gamma$, respectively, and $\mathcal{S}(k)$ is defined by
$$\mathcal{S}(k) = \begin{pmatrix} f_0 \bar{f_0} - 4\Omega^2k^2 & i \text{\upshape Im} f_0 + 2i\Omega k \\ 
i\text{\upshape Im} f_0 - 2i\Omega k & -1 \end{pmatrix}, \qquad k \in \Gamma.$$

\item $\mathcal{M}$ has the asymptotic behavior
$$\mathcal{M}(k) = -\sigma_1 + O(1/k), \qquad k \to \infty.$$
\end{itemize}
\end{proposition}

\begin{figure}
       \vspace{.3cm}
\begin{center}
 \begin{overpic}[width=.4\textwidth]{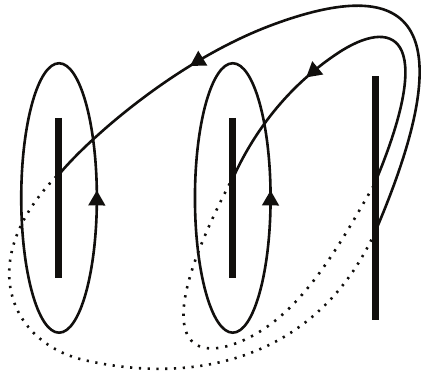}
       \put(25,41){\small $a_1$}
       \put(65.8,41){\small $a_2$}
       \put(40,76){\small $b_1$}
       \put(66,72){\small $b_2$}
       \put(11.5,17.8){\small $k_1$}
      \put(11.5,63.2){\small $\bar{k}_1$}
       \put(52,17.8){\small $k_2$}
      \put(52,63.2){\small $\bar{k}_2$}
     \put(85,72.5){\small $i\bar{z}$}
      \put(81,8.3){\small $-iz$}
        \end{overpic}
     \begin{figuretext}\label{NMcutsright.pdf}
       The homology basis $\{a_j, b_j\}_1^2$ on the Riemann surface $\Sigma_z$ of genus $g = 2$ defined in (\ref{genus2Sigmazdef}).
         \end{figuretext}
     \end{center}
\end{figure}

The auxiliary RH problem presented in Proposition \ref{auxRHprop} can be used to determine the spectral functions $F$ and $G$, which can then be used to set up the main RH problem. In fact, following \cite{MAKNP,NM1993} we can also obtain the solution $f$ directly by combining the main and auxiliary RH problems into a single scalar RH problem. It turns out that the solution can be expressed explicitly in terms of theta functions associated with the Riemann surface $\Sigma_z$ of genus $2$ defined by the equation
\begin{align}\label{genus2Sigmazdef}
y^2 = (k - i\bar{z})(k + iz)\prod_{j = 1}^2 ( k - k_j)(k - \bar{k}_j),
\end{align}
where $k_1, k_2 \in \C$ are such that
$$w^2 + 1 = w_2^2 \prod_{j = 1}^2 (k - k_j)(k - \bar{k}_j), \qquad w(k) := w_2 (k^2 + \rho_0^2),$$
and the solution is parametrized by the two parameters $\rho_0 > 0$ and $w_2 > 0$. 
We view $\Sigma_z$ as a two-sheeted cover of the Riemann $k$-sphere by introducing branch cuts from $k_j$ to $\bar{k}_j$, $j = 1, 2$, and from $-iz$ to $i\bar{z}$; the upper (lower) sheet is characterized by $y/k^3 \to 1$ ($y/k^3 \to -1$) as $k \to \infty$. 
We let $\{a_j, b_j\}_{j = 1}^2$ denote the homology basis on $\Sigma_z$ depicted in Figure \ref{NMcutsright.pdf}. We let $\omega = (\omega_1, \omega_2)^T$, where $\omega_j$, $j =1,2$, is a canonical dual basis of holomorphic one-forms. The associated period matrix $B$ and theta function $\Theta(v |B)$ are defined by
\begin{equation}\label{Thetadef} 
  B_{ij} = \int_{b_j} \omega_i, \quad i,j = 1, 2; \qquad \Theta(v | B) = \sum_{N \in \Z^2} e^{2 \pi i\left(\frac{1}{2} N^T B N + N^T v\right)}, \quad v \in \C^2.
\end{equation}
We denote by $\omega_{PQ}$ the normalized Abelian differential of the third kind on $\Sigma_z$, which has simple poles at the points $P$ and $Q$ with residues $+1$ and $-1$, respectively.

\begin{theorem}\label{diskth}
Let $\rho_0, w_2 > 0$.\footnote{The requirement that the solution be singularity-free imposes further restrictions on these parameters, see \cite{MAKNP}. The correspondence between the parameter $w_2$ used here and the parameter $\mu$ in \cite{MAKNP} is $\mu = w_2 \rho_0^2$.}
Let the function $h(k)$ be defined by
\begin{equation}
  h(k) = \frac{1}{\pi i}\ln\left(\sqrt{w(k)^2 + 1} - w(k)\right), \qquad k \in \Gamma = [-i\rho_0, i\rho_0].
\end{equation}  
Define the $z$-dependent quantities $u \in \C^2$ and $I \in \R$ by
\begin{equation}\label{uIdef}
u = \int_{\Gamma^+} h \omega,
\qquad I = \int_{\Gamma^+} h \omega_{\infty^+\infty^-}.
\end{equation}
Then the function 
\begin{equation}\label{diskernst}
f(z) = \frac{\Theta\left(u - \int_{-iz}^{\infty^-} \omega \bigl | B\right)}{\Theta\left (u + \int_{-iz}^{\infty^-} \omega \bigl | B\right)}e^{I},
\end{equation}
satisfies the BVP denoted by (b) in the introduction with the prescribed value of $\rho_0$.
\end{theorem}

\subsection{A disk rotating around a black hole}
Suppose $f(z)$ is a solution of the BVP formulated in (c) of the introduction and let $\Phi(z,k)$ be the eigenfunction defined by (\ref{ernstlax}) and (\ref{phiinitial}). The boundary values (\ref{holediskconditions}) imply that $\Phi$ has simple poles at the points $k = \pm r_1$ where there regular rotation axis meets the horizon. Thus, the effect of including a central black hole is to add two bound states (which correspond to solitons) to the solution. 

Let $f_0 := f(+i0)$ and $f_1 := f(ir_1)$. The following four propositions are the analogs of Propositions \ref{diskPhionaxisprop}, \ref{eqprop}, \ref{fconstprop}, and \ref{auxRHprop}.

\begin{proposition}\label{FGprop}
The values of $\Phi$ on the $\zeta$-axis can be expressed in terms of two spectral functions $F(k)$ and $G(k)$ as
\begin{align*} 
&\Phi(i\zeta, k^+) =  \begin{pmatrix} \overline{f(i\zeta)} & 1 \\ f(i\zeta) & -1 \end{pmatrix}A(k), & \zeta  > r_1,\quad k \in \hat{\C},
	\\ 
&\Phi(i\zeta, k^+) =  \begin{pmatrix} \overline{f(i\zeta)} & 1 \\ f(i\zeta) & -1 \end{pmatrix}T_1(k)A(k), & 0 < \zeta  < r_1,\quad k \in \hat{\C},
	\\ 
&\Phi(i\zeta, k^+) =  \begin{pmatrix} \overline{f(i\zeta)} & 1 \\ f(i\zeta) & -1 \end{pmatrix}T_2(k)\sigma_1A(k)\sigma_1, & -r_1 < \zeta  < 0,\quad k \in \hat{\C},
	\\ 
&\Phi(i\zeta, k^+) = \begin{pmatrix} \overline{f(i\zeta)} & 1 \\ f(i\zeta) & -1 \end{pmatrix}\sigma_1A(k)\sigma_1, & \zeta  < -r_1,\quad k \in \hat{\C},
\end{align*}
where $A(k)$ is given in (\ref{Adef}) and the $2\times 2$-matrix valued functions $T_1(k)$ and $T_2(k)$ are defined by
\begin{align} \label{T12def}
&T_1(k) = \frac{1}{2 (k - r_1) \Omega_h} \begin{pmatrix} 2(k - r_1) \Omega_h - i f_1 & i \\ - if_1^2 & 2(k - r_1) \Omega_h + i f_1 \end{pmatrix}, \qquad k \in \hat{\C}, 
	\\ \nonumber
&T_2(k) = \overline{T_1(-\bar{k})} , \qquad k \in \hat{\C}.
\end{align}
\end{proposition}

\begin{proposition}\label{holeeqprop}
The spectral functions $F(k)$ and $G(k)$ defined in Proposition \ref{FGprop} satisfy the global relation
\begin{equation}
  \overline{T_1A_+\sigma_1A_+^{-1} T_1^{-1}} = T_2 \sigma_1A_-\sigma_1A_-^{-1}\sigma_1 T_2^{-1}, \qquad k \in \Gamma,
\end{equation}
where $A$ is defined in terms of $F$ and $G$ by equation (\ref{Adef}).
\end{proposition}

\begin{proposition}\label{holefconstprop}
The spectral functions $F(k)$ and $G(k)$ satisfy the relation
\begin{equation}\label{fconstrelation}
  (B^{-1}\Lambda^{-1}\sigma_1\sigma_3\overline{\Lambda B})(\overline{T_1A_+ \sigma_1 A_+^{-1}T_1^{-1}}) = - (T_1 A_+\sigma_1A_+^{-1}T_1^{-1})(B^{-1} \Lambda^{-1}\sigma_1\sigma_3 \overline{\Lambda B}), \quad k \in \Gamma,
\end{equation}
where $B$ is as in (\ref{Bdef}) and $\Lambda(k)$ is defined by
\begin{equation}
 \Lambda(k) = \left(1 - \frac{\Omega}{\Omega_h}\right)I+ \frac{i k \Omega}{\text{\upshape Re}\, f_0}(\sigma_1 - I) \sigma_3.
\end{equation}
\end{proposition}

\begin{proposition}
Suppose $f$ satisfies the BVP denoted by (c) in the introduction. Then the spectral functions $F(k)$ and $G(k)$ are given by
$$F(k) = \mathcal{M}_{12}(k), \qquad G(k) = \mathcal{M}_{22}(k), \qquad k \in \hat{\C},$$
where $\mathcal{M}$ is the unique solution of the following RH problem:
\begin{itemize}
\item \text{$\mathcal{M}(k)$ is analytic for $k \in \hat{\C} \setminus (\Gamma \cup \{-r_1, r_1\})$.} 

\item  Across $\Gamma$, $\mathcal{M}(k)$ satisfies the jump condition (\ref{calMjump}), where $\mathcal{S}(k)$ is defined by
\begin{equation}
\mathcal{S}(k) = T_1^{-1}B^{-1}\Lambda^{-1}\sigma_1\sigma_3\overline{\Lambda B} T_2\sigma_1, \qquad k \in \Gamma.
\end{equation}

\item $\mathcal{M}$ has the asymptotic behavior $\mathcal{M}(k) = \sigma_1 + O(k^{-1})$ as $k \to \infty$.

\item The entries of $\mathcal{M}$ have simple poles at $k = r_1$ and $k = -r_1$. The associated residues are given by
\begin{equation}\label{Mresidues}
\underset{r_1}{\text{\upshape Res}} \, \mathcal{M}(k) = \frac{1}{\alpha}\begin{pmatrix} -f_1 & 1 \\ -f_1^2 & f_1 \end{pmatrix}, \qquad \underset{-r_1}{\text{\upshape Res}} \, \mathcal{M}(k) = \frac{1}{\alpha}\begin{pmatrix} f_1 & -|f_1|^2 \\ f_1/\bar{f}_1 & -f_1 \end{pmatrix},
\end{equation}
where 
\begin{equation}\label{alphadef}  
  \alpha = \frac{d^+}{d\zeta}\biggl|_{\zeta = r_1} \re f(i\zeta)
\end{equation}
and $d^+/d\zeta$ denotes the right-sided derivative.
\end{itemize}
\end{proposition}

It turns out that the main and auxiliary matrix RH problems can be combined into a single scalar RH problem formulated on the Riemann surface $\Sigma_z$ of genus $4$ consisting of all points $(k, y) \in \C^2$ such that
\begin{equation}\label{genus4Sigmazdef}  
  y^2 = (k + iz) (k - i\bar{z}) \prod_{j = 1}^4 (k - k_j) (k - \bar{k}_j),
\end{equation}
where $k_1, \dots, k_4 \in \C$ are such that
$$w^2 + 1 = \frac{w_4^2 \prod_{j = 1}^4 (k - k_j)(k - \bar{k}_j)}{(k^2 - r_1^2)^2}, \qquad
w(k) := \frac{w_4 k^4 + w_2 k^2 + \rho_0^2(w_2 - w_4 \rho_0^2)}{(k^2 - r_1^2)},$$
and the solution is parametrized by the four real parameters $\rho_0$, $r_1$, $w_2$, and $w_4$.

\begin{figure}
       \vspace{.3cm}
\begin{center}
 \begin{overpic}[width=.5\textwidth]{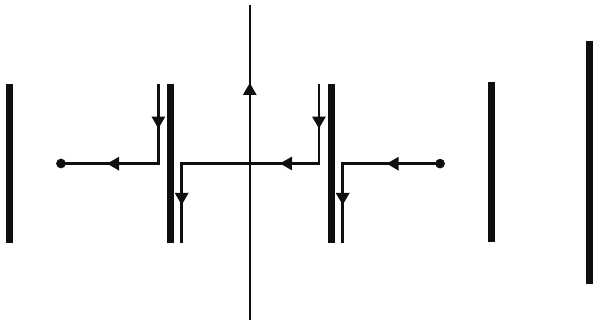}
        \put(0,9.2){\small $k_1$}
      \put(0,42){\small $\bar{k}_1$}
       \put(26.5,9.2){\small $k_2$}
      \put(26.5,42){\small $\bar{k}_2$}
       \put(53.5,9.2){\small $k_3$}
      \put(53.5,42){\small $\bar{k}_3$}
       \put(80.5,9.2){\small $k_4$}
      \put(80.5,42){\small $\bar{k}_4$}
     \put(96.3, 49){\small $i\bar{z}$}
      \put(93,2.5){\small $-iz$}
        \put(43,40){\small $\Gamma^+$}
        \put(46.5,23.3){\small $\gamma$}
       \end{overpic}
     \begin{figuretext}\label{contours1.pdf}
       The genus $4$ Riemann surface $\Sigma_z$ defined in (\ref{genus4Sigmazdef}) presented as a two-sheeted cover of the complex $k$-plane together with the contours $\gamma$ and $\Gamma^+$.
         \end{figuretext}
     \end{center}
\end{figure}

We let $\gamma$ denote the contour on $\Sigma_z$ which projects to the contour\footnote{For $n$ complex numbers $\{a_j\}_1^n$, we let $[a_1, \dots, a_n]$ denote the directed contour $\cup_{j = 1}^{n-1} [a_j, a_{j+1}]$.}
\begin{equation}\label{gammaprojection}
[r_1, \re k_3 + \epsilon, k_3 + \epsilon]
\cup [\bar{k}_3 - \epsilon, \re k_3 - \epsilon, \re k_2 + \epsilon, k_2 + \epsilon]
\cup [\bar{k}_2 - \epsilon, \re  k_2 - \epsilon, -r_1] 
\end{equation}
in the complex $k$-plane, where $\epsilon > 0$ is an infinitesimally small positive number, and which lies in the upper sheet for $\re  k < \zeta$ and in the lower sheet for $\re k > \zeta$, see Figure \ref{contours1.pdf}.
For simplicity, we assume that
\begin{equation}\label{r1betweencuts}
  0 < \re  k_3 < r_1 < \re  k_4.
\end{equation}
Using notations analogous to those used in Theorem \ref{diskth}, we can state the following result.

\begin{theorem}{\upshape \bf \cite{L}}\label{mainth}
Let $\rho_0, r_1, w_2, w_4$ be strictly positive numbers such that (\ref{r1betweencuts}) holds; the requirement that the solution be singularity-free imposes further restrictions on these parameters, see \cite{L}.
Let the function $h(k)$ be defined by
\begin{equation}\label{hdef}
  h(k) = \frac{1}{\pi i}\ln\left(\sqrt{w(k)^2 + 1} - w(k)\right), \qquad k \in \Gamma = [-i\rho_0, i\rho_0].
\end{equation}  
Define the $z$-dependent quantities $u \in \C^4$ and $I \in \R$ by
$$u = \int_{\Gamma^+} h \omega + \int_\gamma \omega,
\qquad I = \int_{\Gamma^+} h \omega_{\infty^+\infty^-} + \int_\gamma \omega_{\infty^+\infty^-}.$$
Then the function 
\begin{equation}\label{ernstsolution}
f(z) = \frac{\Theta\left(u - \int_{-iz}^{\infty^-} \omega \bigl | B\right)}{\Theta\left (u + \int_{-iz}^{\infty^-} \omega \bigl | B\right)}e^{I},
\end{equation}
satisfies the BVP denoted by (c) in the introduction with the prescribed values of $\rho_0$ and $r_1$, and with the values of $\Omega_h$, $\Omega$, and $\text{\upshape Re} f_\Omega(+i0)$ given by
\begin{align*}
\Omega_h =  -\frac{1}{a_{hor}}, \quad 
\Omega =  \frac{w_4 \Omega_h \text{\upshape Re}f_0 +  \sqrt{-2 w_4 \Omega_h^4 \text{\upshape Re}f_0}}{w_4\text{\upshape Re}f_0+2 \Omega_h^2}, \quad  	
\text{\upshape Re}f_\Omega(+i0) =  \text{\upshape Re}(f_0)\Bigl(1 - \frac{\Omega}{\Omega_h}\Bigr)^2,
\end{align*}
where $a_{hor}\in \R$ denotes the (necessarily constant) value of the metric function $a$ on the horizon. 
\end{theorem}

\begin{remark}
Explicit expressions in terms of theta functions for $a_{hor}$ and $f_0$ as well as for the metric functions $e^{2U}, a, e^{2\kappa}$ corresponding to the Ernst potential (\ref{ernstsolution}) are presented in \cite{L}.
\end{remark}

\subsection{A disk with a Neumann condition}
Suppose $f(z)$ is a solution of the BVP formulated in (d) of the introduction for some choice of the parameters $\rho_0 > 0$ and $\Omega > 0$ such that $2\Omega \rho_0 < 1$.
Then the following analogs of Propositions \ref{fconstprop} and \ref{auxRHprop} are valid (see \cite{LP2019} for details). 

\begin{proposition}\label{Neumannfconstprop}
Suppose that $\partial_\zeta f_\Omega = 0$ on the disk. Then the spectral functions $F(k)$ and $G(k)$ satisfy the relation
\begin{equation*}
  (B^{-1}\Lambda^{-1}\sigma_1\overline{\Lambda B})(\overline{A_+ \sigma_1 A_+^{-1}}) = - (A_+\sigma_1A_+^{-1})(B^{-1} \Lambda^{-1}\sigma_1 \overline{\Lambda B}), \quad k \in \Gamma,
\end{equation*}
where $B$ and $\Lambda(k)$ are given in (\ref{Bdef}).
\end{proposition}

\begin{proposition}\label{NeumannauxRHprop}
Suppose $f$ satisfies the BVP denoted by (d) in the introduction. 
Then the spectral functions $F(k)$ and $G(k)$ are given by
$$F(k) = \mathcal{M}_{12}(k), \qquad G(k) = \mathcal{M}_{22}(k), \qquad k \in \C,$$
where $\mathcal{M}$ is the unique solution of the following RH problem:
\begin{itemize}
\item \text{$\mathcal{M}(k)$ is analytic for $k \in \C \backslash \Gamma$.} 

\item  Across $\Gamma$, $\mathcal{M}(k)$ satisfies the jump condition (\ref{calMjump}), where $\mathcal{S}(k)$ is defined by
$$\mathcal{S}(k) = \begin{pmatrix} 0 & 1 \\ -1 & 4i\Omega k \end{pmatrix}, \qquad k \in \Gamma.$$

\item $\mathcal{M}$ has at most logarithmic singularities at the endpoints of $\Gamma$.

\item $\mathcal{M}$ has the asymptotic behavior $\mathcal{M}(k) = \sigma_1 + O(k^{-1})$ as $k \to \infty$.
\end{itemize}
\end{proposition}

\begin{figure}
       \vspace{.3cm}
\begin{center}
 \begin{overpic}[width=.35\textwidth]{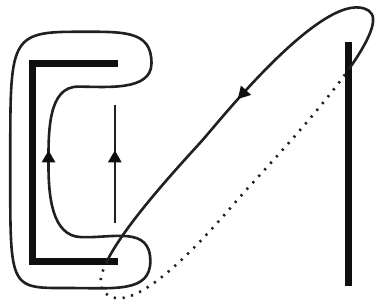}
       \put(33,37.5){\small $\Gamma^+$}
       \put(15.5,38.5){\small $a$}
       \put(59,56){\small $b$}
       \put(32.5,10.5){\small $k_1$}
      \put(32.5,63){\small $\bar{k}_1$}
     \put(89,72){\small $i\bar{z}$}
      \put(84.5,-.6){\small $-iz$}
        \end{overpic}
     \begin{figuretext}\label{Neumanncuts.pdf}
       The homology basis $\{a, b\}$ on the Riemann surface $\Sigma_z$ of genus $g = 1$ defined in (\ref{NeumannSigmazdef}) and the contour $\Gamma^+$.
         \end{figuretext}
     \end{center}
\end{figure}

The solution of the auxiliary RH problem of Proposition \ref{auxRHprop} yields the following explicit formulas for $F$ and $G$, which can be used to set up the main RH problem:
\begin{align*}
\begin{cases}
F(k) 
= \frac{i}{4\Omega \mu(k)}\Big(\frac{d_1(k)}{E(k)} - \frac{E(k)}{d_1(k)}\Big),
	\vspace{1mm} \\ 
G(k)= \frac{i}{4\Omega\mu(k)}\Big( \frac{1}{E(k)}- E(k)\Big),
\end{cases} \quad
k \in \C \setminus \Gamma,
\end{align*}
where $\mu(k) = \sqrt{k^2 + 1/(4\Omega)^2}$, $d_1(k) = 2\Omega i(k - \mu(k))$, and 
\begin{align}\label{Ehdef}
E(k) = e^{\mu(k) \int_{\Gamma} \frac{h(s)}{\mu(s)} \frac{ds}{s - k}} \quad \text{with} \quad h(k) = -\frac{\text{\upshape arcsin}(2 i \Omega k)}{\pi}.
\end{align}
Furthermore, by combining the main and auxiliary RH problems into a single scalar RH problem, the solution $f$ can be expressed in terms of theta functions associated with the family of Riemann surfaces $\Sigma_z$ of genus $1$ defined by 
\begin{align}\label{NeumannSigmazdef}
y^2 = (k - i\bar{z})(k + iz)( k - k_1)(k - \bar{k}_1),
\end{align}
where $k_1 = -i/(2\Omega)$.
We view $\Sigma_z$ as a two-sheeted cover of the Riemann $k$-sphere by introducing a vertical branch cut from $-iz$ to $i\bar{z}$ and a branch cut from $k_1$ to $\bar{k}_1$ which passes to the left of $\Gamma$, see Figure \ref{Neumanncuts.pdf}; the upper (lower) sheet is characterized by $y/k^2 \to 1$ ($y/k^2 \to -1$) as $k \to \infty$. 
Let $\{a, b\}$ denote the homology basis on $\Sigma_z$ shown in Figure \ref{Neumanncuts.pdf}. Let $\omega$ denote the unique holomorphic one-form on $\Sigma_{z}$ such that $\int_{a}\omega=1$ and let $B := \int_{b}\omega \in \C$.

\begin{theorem}{\upshape \bf \cite{LP2019}} \label{Neumanndiskth}
Let $\rho_0 > 0$ and $\Omega > 0$ be such that $2\Omega \rho_0 < 1$.
Define the $z$-dependent quantities $u \in \C^2$ and $I \in \R$ by (\ref{uIdef}), where the integrals are contour integrals on the genus $1$ surface  $\Sigma_z$ defined in (\ref{NeumannSigmazdef}) and $h(k)$ is given by (\ref{Ehdef}).
Then the function $f(z)$ defined by (\ref{diskernst}) satisfies the BVP denoted by (d) in the introduction with the prescribed values of $\rho_0$ and $\Omega$.
\end{theorem}

 \bigskip
\noindent
{\bf Acknowledgement} {\it ASF acknowledges support from the Guggenheim foundation and from the EPSRC in the form of a senior fellowship. JL acknowledges support from the G\"oran Gustafsson Foundation, the Ruth and Nils-Erik Stenb\"ack Foundation, the Swedish Research Council, Grant No. 2015-05430, and the European Research Council, Grant Agreement No. 682537.}

\bibliography{is}

\end{document}